\documentclass[10pt]{article}
\usepackage{amsmath,amssymb,bbm,times}

\makeatletter
\newfont{\footsc}{cmcsc10 at 8truept}
\newfont{\footbf}{cmbx10 at 8truept}
\newfont{\footrm}{cmr10 at 10truept}

\newtheorem{Definition}{Definition}

\newtheorem{Lemma}{Lemma}

\newtheorem{Theorem}{Theorem}

\newtheorem{Conjecture}{Conjecture}
\newenvironment{proof}{
  \noindent\textbf{Proof}\ }{\hspace*{\fill}
  \begin{math}\Box\end{math}\medskip}
\begin{document}
\title{\Large\bf Explicit tough Ramsey graphs}

\author{\begin{tabular}[t]{c@{\extracolsep{8em}}c} 
  Le Anh Vinh & Dang Phuong Dung\\
  \textit{Mathematics Department}      & \textit{School of Banking and Finance} \\
  \textit{Harvard University}     & \textit{University of New South Wales} \\
  \textit{Cambridge, MA 02138}               &  \textit{Sydney 2052 NSW, Australia}\\
  \textit{vinh@math.harvard.edu}         & \textit{phuogdug@yahoo.com}
\end{tabular}}

\date{\empty}
\maketitle

\begin{center} {\em Dedicated to our parents on the occasion of their 30th anniversary.} \end{center}

\section*{\centering \large Abstract}

{\em
A graph $G$ is $t$-\textit{tough} if any induced subgraph of it with $x > 1$ connected components is obtained from $G$ by deleting at least $tx$ vertices. Chv\'atal conjectured that there exists an absolute constant $t_0$ so that every $t_0$-tough graph is pancyclic. This conjecture was disproved by Bauer, van den Heuvel and Schmeichel by constructing a $t_0$-tough triangle-free graph for every real $t_0$. 
For each finite field $\mathbbm{F}_q$ with $q$ odd, we consider graphs associated to the finite Euclidean plane and the finite upper half plane over $\mathbbm{F}_q$. These graphs have received serious attention as they have been shown to be Ramanujan (or asymptotically Ramanujan) for large $q$. We will show that for infinitely many $q$, these graphs provide further counterexamples to Chv\'atal's conjecture. They also provide a good constructive lower bound for the Ramsey number $R(3,k)$.
}

\begin{center}
\small Mathematics Subject Classifications: 05C35, 05C38, 05C55, 05C25.

\small Keywords: finite Euclidean graph, finite non-Euclidean graph, Ramsey graphs,\\ tough graphs, girth, triangle-free, diameter.
\end{center}

\section*{\large 1. Introduction}
\label{sec:in}

The \textit{toughness} $t(G)$ of a graph $G$ is the largest real $t$ so that for every positive integer $x \geq 2$ one must delete at least $tx$ vertices from $G$ in order to get an induced subgraph of it with at least $x$ connected components. $G$ is $t$-tough if $t(G) \geq t$. This parameter was introduced by Chv\'atal in \cite{chvatal}. Chv\'atal proposed the following conjecture.

\begin{Conjecture} (\cite{chvatal}) A simple unlabeled graph on $n$ vertices is called \textit{pancyclic} if it contains cycles of all lengths $3, 4,\ldots, n$. Then there exists an absolute constant $t_0$ such that every $t_0$-tough graph is pancyclic. 
\end{Conjecture}

This conjecture was disproved by Bauer, van den Heuvel and Schmeichel \cite{bauer} who constructed, for every real $t_0$, a $t_0$-tough triangle-free graph. They defined a sequence of triangle-free graphs $H_1,H_2,H_3,\ldots$ with $|V(H_j)| = 2^{2j-1} (j+1)!$ and $t(H_j) \geq \sqrt{2j+4}/2$. We present here two new explicit constructions based on finite analogues of Euclidean and non-Euclidean planes. 

The first construction is based on a remarkable new approach of Wildberger to trigonometry and Euclidean geometry. In \cite{quadrance}, Wildberger replaces distance by quadrance and angle by spread, thus allowing the development of Euclidean geometry over any field. The following definition follows from \cite{quadrance}.

\begin{Definition}
  Let $\mathbbm{F}_q$ be a finite field of order $q$ where $q$ is an odd prime power. The quadrance $Q ( X, Y )$ between the points $X = ( x_1,
  x_2)$, and $Y = ( y_1, y_2 )$ in $\mathbbm{F}_q^2$ is the number
  \begin{equation} Q ( X, Y ) = ( y_1 - x_1 )^2 + ( y_2 - x_2 )^2 . \end{equation}
\end{Definition}

Let $q$ be an odd prime power and $\mathbbm{F}_q$ be the finite field with $q$ elements. We can associate graphs to the symmetric plane $\mathbbm{F}_q^2$ using the quadrance. 
\begin{Definition}
For a fixed $a \in \mathbbm{F}_q$, the \textit{quadrance graph} $D_q(a)$ has the vertex set $\mathbbm{F}_q^2$, and $X, Y \in \mathbbm{F}_q^2$ are adjacent if and only if $Q ( X, Y ) = a$. 
\end{Definition}

This graph or the so-called finite Euclidean graph was first studied by Moorhouse \cite{moorhouse} and then by Medrano et al \cite{before}. Note that for any $a, b \neq 0 \in \mathbbm{F}_q$ then $D_q(a)$ and $D_q(b)$ are isomorphic. In usual $2$-dimensional Euclidean space $\mathbbm{R}^2$, then the quadrance between $X, Y$ is unity if and only if the distance between $X, Y$ is unity, so the unit-quadrance graph $D_q(1)$ becomes the unit-distance graph. 

The second construction is based on the well-known finite upper half planes constructed in a similar way using an analogue of Poincar\'e's non-Euclidean distance. We follow the construction in \cite{survey}. Let $\mathbbm{F}_q$ be the finite field with $q = p^r$ elements, where $p$ is an odd prime. Suppose $\sigma$ is a generator of the multiplicative group $\mathbbm{F}_q^*$ of nonzero elements in $\mathbbm{F}_q$. The extension $\mathbbm{F}_q(\sigma)$ is analogous to $\mathbbm{C}=\mathbbm{R}[i]$. We define the \textit{finite Poincar\'e upper half-plane} as 
\begin{equation} H_q = \{z=x+y\sqrt{\sigma}: x, y \in \mathbbm{F}_q \, \text{and} \, y \neq 0 \}.\end{equation}
Note that ``half-plane'' is something of a misnomer since $y \neq 0$ may not be a good finite analogue of the condition $y > 0$ that defines the usual Poincar\'e upper half-plane in $\mathbbm{C}$. In fact, $H_q$ is more like a double covering of a finite upper half-plane. We use the familiar notation from complex analysis for $z = x+y \sqrt{\sigma} \in H_q$: $x= Re(z)$, $y= Im(z)$, $\bar{z} = x- y\sqrt{\sigma} = z^q$, $N(z) = $ Norm of $z = z\bar{z} = z^{1+q}$. The \textit{Poincar\'e distance} between $z, w \in H_q$ is
\begin{equation} d(z,w) = \frac{N(z-w)}{Im(z)Im(w)}.
\end{equation}
This distance is not a metric in the sense of analysis, but it is $GL(2,\mathbbm{F}_q)$-invariant: $d(gz,gw) = d(z,w)$ for all $g \in GL(2,\mathbbm{F}_q)$ and all $z,w \in H_q$. We can attach graphs to $H_q$ by a method analogous to that which led to the unit-quadrance graphs $D_q$. 

\begin{Definition}
For a fixed $a \in \mathbbm{F}_q$, the \textit{finite non-Euclidean graph} $V_q(\sigma,a)$ has vertices as the points in $H_q$ and edges between vertices $z,w$ if and only if $d(z,w) = a$. 
\end{Definition}

Except when $a=0$ or $a=4\sigma$, $P_q(\sigma,a)$ is a connected $(q+1)$-regular graph. When $a=0, 4\sigma$ then $P_q(\sigma,a)$ is disconnected, with one or two nodes, respectively, per connected component. As $a$ varies, we have $q-2$ $(q+1)$-regular graphs $P_q(\sigma,a)$. The question of whether these graphs are always nonisomorphic or not is still open. These graphs are studied extensively in the literature, see [1-2, 8, 9, 11, 12, 20-22].

The rest of this paper is organized as follows. In Section 2 we determine the girth, diameter and triangle-freeness of these graphs. In Section 3 we establish some useful
facts about toughness and cutsets of graphs that we will need in the study of quadrance graphs and finite non-Euclidean graphs. We then show in Sections 4 and 5 that for infinitely many values of $q$, these graphs constitute counterexamples to Chv\'atal's conjecture and also provide a good constructive lower bound for the Ramsey number $R(3,k)$. Finally, we make some further remarks on the chromatic numbers of these graphs.

\section*{\large 2. Girth and diameter}

We define the girth to be the length of a shortest circuit in a graph. The diameter is the maximum length of shortest paths between two vertices of a graph. The main aim of this section is to determine the girth and diameter of finite Euclidean and non-Euclidean graphs.

We first consider the finite Euclidean graph $D_q$.

\begin{Definition}
The circle $C_k(A_0)$ in $\mathbbm{F}_q^2$ with center $A_0 \in \mathbbm{F}_q \times \mathbbm{F}_q$ and quadrance $k$ is the set of all points $X \in \mathbbm{F}_q \times \mathbbm{F}_q$ such that
\[Q(A_0,X)=k.\]
\end{Definition}

The following lemma (\cite{vinh4, quadrance}) gives us the number of intersections between any two circles in $\mathbbm{F}_q^2$.

\begin{Lemma} \label{khac khong}
Let $i, j \ne 0$ in $\mathbbm{F}_q$ and let $X, Y$ be two distinct points in $\mathbbm{F}_q^2$ such that $Q(X,Y) = k \ne 0 $. Then $\left|C_i(X) \cap C_j(Y)\right|$ only depends on $i, j$ and $k$. Precisely,  let \[f ( i, j, k ) = i j - (k - i - j )^2/4.\] Then the number of intersection points is $p_{i j}^k$, where
  \begin{equation}\label{eq2}
  p_{ij}^k =
   \begin{cases}
    0  & \text{if}\; \; f(i,j,k)\; \text{is a non-square},\\
    1  & \text{if}\; \; f(i,j,k) = 0,\\
    2 &\text{if}\; \; f(i,j,k)\; \text{is a square}.
   \end{cases}
  \end{equation}

\end{Lemma}

If $-1$ is a square in $\mathbbm{F}_q$ then there exist distinct points $X, Y$ in $\mathbbm{F}_q^2$ with $Q(X,Y)=0$. The following lemma gives us the intersection number of two circles $C_i(X)$ and $C_j(Y)$ in this case.

\begin{Lemma}\label{bang khong} (\cite{vinh4})
For any $i,j \neq 0$ in $\mathbbm{F}_q$. Suppose that $X, Y$ are two distinct points in $\mathbbm{F}_q^2$ with $Q(X,Y)=0$. Then the circle $C_i(X)$ intersects $C_j(Y)$ if and only if $i \neq j$. Furthermore, if $i \neq j$ then two circles intersect at only one point.
\end{Lemma}

The proof of Lemma \ref{bang khong} is similar to the proof of Lemma \ref{khac khong}.

\begin{Theorem}\label{girth}
For a fixed $a \neq 0 \in \mathbbm{F}_q$, the girth of quadrance graph $Q_q(a)$ is $3$ if $3$ is a square in $\mathbbm{F}_q$ and $4$ otherwise.
\end{Theorem}

\begin{proof} We have $\{(0,0), (0,1),(1,1),(1,0)\}$ is a cycle of length $4$ in $Q_q(1)$ so the girth of unit-quadrance graph is at most $4$. Besides, $Q_q(a), Q_q(b)$ are isomorphic for any $a, b \neq 0 \in \mathbbm{F}_q$, so all quadrance graphs $Q_q(a)$ have girth at most $4$. From Lemma \ref{khac khong}, the girth of quadrance graph $Q_q(a)$ is $3$ if and only if $f(a,a,a) = 3a^2/4$ is a square in $\mathbbm{F}_q$ or $3$ is a square in $\mathbbm{F}_q$. This concludes the proof of the theorem. 
\end{proof}

\begin{Theorem}\label{diameter}
If $q \equiv 3$ mod $4$ then $D_q$ has diameter $3$. Otherwise, $D_q(a)$ has diameter $3$ or $4$.
\end{Theorem}

\begin{proof} Since all quadrance graphs are isomorphic, we only need to prove the theorem for unit-quadrance graph $D_q = D_q(1)$. 
Define $g(x) = 4f(x,1,1) = (4-x)x$. Then there exist $u, v \neq 0 \in \mathbbm{F}_q$ such that $g(u)$ is a square in $\mathbbm{F}_q$ while $g(v)$ is not. Choose any $X \in C_v((0,0))$. By Lemma \ref{khac khong}, $X$ is not reachable from $(0,0)$ by two unit-steps. Thus, the diameter of $D_q$ is at least $3$. 

Now suppose that $q \equiv 3$ mod $4$. Then $D_q$ is a $(q+1)$-regular graph. For any vertex $X \in D_q$, we define the set $V(X)$ to be all vertices that are adjacent to $X$ and $T(X)$ to be all the vertices that are connected by a path of at most two edges from $X$. Then 
\begin{equation}
T(X) = \bigcup_{Y \in V(X)} V(Y).
\end{equation}
By Lemma \ref{khac khong}, a fixed vertex $W$ can be in $V(Y)$ for at most $2$ vertices $Y \in V(X)$. Thus, $|T(X)| \geq (q+1)q/2$. This implies that $T(X)$ contains elements from no fewer than $\left\lceil q/2 \right\rceil$ circles $C_i(X)$. Thus, there are no more than $q - (q+1)/2 = (q-1)/2$ circles outside of $T(X)$. Let $Z$ be any point then $Z$ is adjacent to at most $q-1$ points that are not in $T(X)$. Thus $Z$ is connected to $X$ by a path of length at most three. This implies that $D_q$ has diameter $3$. 

Suppose that $q \equiv 1$ mod $4$. For any two points $X, Y$ with $Q(X,Y) \neq 1$, by Lemma \ref{bang khong} we can choose a point $Z$ such that $Q(X,Z)=1, Q(Y,Z)=0$. By Lemma \ref{bang khong} again, we can choose a point $W$ such that $Q(Z,W) = 1$ and $Q(W,Y) = u$. Since $g(u)$ is a square in $\mathbbm{F}_q$, by Lemma \ref{khac khong}, we can choose a point $T$ such that $Q(W,T)=Q(T,Y)=1$. Thus, $Y$ is connected to $X$ by a path of length at most four. If $Q(X,Y) = 1$ then $Y$ is adjacent to $X$ in the graph. This concludes the proof of the theorem.
\end{proof}

The girth and diameter of finite non-Euclidean graphs have been studied before. The following theorem is due to Celniker \cite{celniker}.

\begin{Theorem} (\cite{celniker}) \label{girth1} Let $q$ be odd and $a \notin \{0,4\sigma\}$. The girth of of $V_q(\sigma,a)$ is either $3$ or $4$. Furthermore, the girth is $3$ if $a=2\sigma$ and $q \equiv 3$ mod $4$ or if $a$ and $a-3\sigma$ are squares in $\mathbbm{F}_q$. The girth is $4$ if $a=2\sigma$ and $q \equiv 1$ mod $4$.
\end{Theorem}

In \cite{angel}, Angel and Evans obtained the following result for the diameter of finite non-Euclidean graphs $V_q(\sigma,a)$.

\begin{Theorem} (\cite{angel}) \label{diameter1} Let $q$ be odd and $a \notin \{0,2\sigma, 4\sigma\}$. Then $V_q(\sigma,a)$ has diameter $3$ or $4$ according to whether $\sigma -a$ is a square or a non-square in $\mathbbm{F}_q$. Besides, the diameter of $X_{2\sigma}$ is $3$, unless $q =3$ or $5$, in which case the diameter is $2$. 
\end{Theorem}

\section*{\large 3. Toughness and cutsets}

Let $G$ be a $d$-regular graph with $n$ eigenvalues $\lambda_1 \geq \lambda_2 \geq \ldots \geq \lambda_n$. It is well-known that the largest eigenvalue $\lambda_1 = d$. We know that for $a \neq 0$ the quadrance graph $D_q(a)$ is a regular graph with degree $\Delta(D_q(a)) = q-(-1)^{(q-1)/2}$ (see \cite{before}). Medrano et al. \cite{before} give a general bound for eigenvalue of $D_q(a)$.

\begin{Lemma}\label{eigenvalue} \cite{before} Let $\lambda \ne \Delta(D_q(a))$ be any eigenvalue of graph $D_q(a)$. Then 
\begin{equation}\label{eigen in} | \lambda | \leq 2 q^{1/2}. \end{equation}
\end{Lemma}

The finite non-Euclidean graphs $V_q(\sigma,a)$ are also known as $(q+1)$-regular graphs for $a \neq 0, 4\sigma$. By combining works of Weil \cite{weil}, Evans \cite{evans1, evans2}, Katz \cite{katz1, katz2}, Li \cite{li} and many others, it was proved that the finite non-Euclidean graphs are all Ramanujan graphs, that is, we have the following bound for eigenvalue of $V_q(\sigma,a)$.

\begin{Lemma}\label{eigenvalue 2} \cite{survey} Let $\lambda \ne q+1 $ be any eigenvalue of graph $V_q(\sigma,a)$. Then 
\begin{equation} \label{eigen in 2} | \lambda | \leq 2 q^{1/2}. \end{equation}
\end{Lemma}

An $(n,d,\mu)$-\textit{graph} is a $d$-regular graph on $n$ vertices, in which every eigenvalue $\mu < d$ satisfies $|\mu| \leq \lambda$. The following result is due to Alon in \cite{alon1}.

\begin{Theorem} \cite{alon1} \label{tough theorem} Let $G = (V,E)$ be an $(n,d,\lambda)$-graph. Then the toughness $t=t(G)$ of $G$ satisfies
\begin{equation} \label{tough in}
t > \frac{1}{3} ( \frac{d^2}{\lambda d + \lambda^2} - 1 ).
\end{equation}
\end{Theorem}

Let $\sigma$ be an arbitrary orientation of graph $G$, and let $D$ be the incidence matrix of $G^{\sigma}$. Then the \textit{Laplacian} of $G$ is the matrix $Q(G) = DD^{T}$. It is easy to show that the Laplacian does not depend on the orientation $\sigma$, and hence is well-defined. Let $G$ be a $d$-regular graph. If the adjacency matrix $A$ of $G$ has eigenvalues $\lambda_1 \geq \ldots \geq \lambda_n$, then the Laplacian $Q$ has eigenvalues $\theta_1 = d-\lambda_1 \leq \ldots \leq \theta_n = d-\lambda_n$. It is well-known (see \cite{algebraic graph theory}, pages 287-288) that, if $G$ is a graph on $n$ vertices then the minimum value of 
\begin{equation}\label{ol}
\frac{\sum_{uv \in E(G)} (x_u-x_v)^2}{\sum_u x_u^2},
\end{equation}
as $x$ ranges over all nonzero vectors orthogonal to $(1,\ldots,1)$, is $\theta_2(G)$. The maximum value is $\theta_n(G)$. 

If $S \subset V(G)$, let $\delta S$ denote the set of edges with one end in $S$ and the other in $V(X)-S$. Suppose that $|S| = s$. Let $z$ be the vector whose value is $n-s$ on the vertices in $S$ and $-s$ on the vertices not in $S$. Then $z$ is orthogonal to $(1,\ldots,1)$, so we have
\begin{equation}\label{bis}
\theta_2(G) \leq \frac{\sum_{uv \in E(G)} (x_u-x_v)^2}{\sum_u x_u^2} = \frac{|\delta S|n^2}{s(n-s)^2+(n-s)s^2}=\frac{n|\delta S|}{|S|(n-|S|)} \leq \theta_n(G).
\end{equation}

The \textit{bisection width} of graph $G$ on $n$ vertices is the minimum value of $|\delta S|$, for any subset $S$ of size $\lfloor n/2 \rfloor$. Let $bip(G)$ denote the maximum number of edges in a spanning bipartite subgraph of $G$. This equals the maximum value of $|\delta S|$, where $S$ ranges over all subsets of $V(X)$ with size at most $|V(X)/2$. We have the following lemma which is an immediate consequence of (\ref{bis})

\begin{Lemma}\label{bip} (\cite{algebraic graph theory}) If $G$ is a graph with $n$ vertices, then
\begin{enumerate}
	\item the bisection width of $G$ is at least $n \theta_2(G) / 4 (1+o(1))$, and
	\item $bip(X) \leq n \theta_n(G)/4$. 
\end{enumerate}

\end{Lemma}

\section*{\large 4. Properties of quadrance graphs}

The following theorem summarizes some of the properties of quadrance graphs. 

\begin{Theorem} \label{main theorem} Let $q$ be any prime of the form $q = 12k + 7$ with $k \geq 0$ and $a \neq 0$. The quadrance graph $D_q(a)$ is a $d_q=(q + 1)$-regular graph on $n_q=q^2$ vertices with the following properties. 
\begin{enumerate}
	\item $D_q(a)$ is triangle-free and has diameter $3$.
	\item The toughness of $D_q(a)$ is at least $q^{1/2}/6 = n_q^{1/4}/6$.
	\item The bisection width of $D_q(a)$ is at least $q^2 (q-2q^{1/2})/(4+o(1))$.
	\item The maximum number of edges in a spanning bipartite subgraph of $D_q(a)$ is at most $bip(D_q(a)) \leq q^2 (q+2q^{1/2}) /4$.
	\item The independence number of $D_q(a)$ is at most $2q^{3/2} = 2n_q^{3/4}$, and hence its chromatic number is at least $n_q^{1/4}/2$ 
\end{enumerate}
\end{Theorem}

\begin{proof}
From Lemma \ref{khac khong}, $D_q(a)$ has no triangle as $f(a,a,a) = 3a^2/4$ is a non-square in $\mathbbm{F}_q$ (If $q = 12k+7$ then $3$ is a non-square in $\mathbbm{F}_q$). Part 1 now follows from Theorem \ref{diameter}. 

Part 2 follows directly from Lemma \ref{eigenvalue} and Theorem \ref{tough theorem}. Parts 3 and 4 are immediate from Lemma \ref{eigenvalue} and Lemma \ref{bip}.

Part 5 follows easily from (\ref{bis}) as follows. Let $S$ be a maximum independent set of $D_q(a)$, then $|\delta S| = |S| (q+1)$. From (\ref{bis}), we have
\[ \frac{q^2 |S| (q+1)}{|S|(q^2-|S|)} \leq \theta_n(D_q(a)) \leq q+1+2q^{1/2}. \]
This implies that the independence number $\alpha(D_q(a))$ is at most 
\begin{equation}\label{idd bound}
\alpha(D_q(a))=|S| \leq q^2 - \frac{q^2(q+1)}{q+1+2q^{1/2}} \leq 2 q^{3/2} = 2 n_q^{3/4}.
\end{equation}
Hence, the chromatic number $\chi(D_q(a))$ is at least
\begin{equation}\label{chromatic number bound}
\chi(D_q(a)) \geq \frac{|V(D_q(a))|}{\alpha(D_q(a))} \geq \frac{n_q}{2 n_q^{3/4}} = n_q^{1/4}/2.
\end{equation}
This concludes the proof of the theorem.
\end{proof}

Theorem \ref{main theorem} shows that the quadrance graph $D_q(a)$, where $q$ is a prime of the form $q=12k+7$ and $a \neq 0 \in \mathbbm{F}_q$, is an explicit triangle-free graph on $n_q=q^2$ vertices whose chromatic number exceeds $0.5n_q^{1/4}$. Note that the lower bound was already observed in \cite{old}. In addition, the quadrance graph $D_q(a)$ is an explicit construction showing that $R(3,k) \geq \Omega(k^{4/3})$. Moreover, $D_q(a)$ has $q^2(q+1)/2$ edges so the graph $D_q(a)$ is also an explicit construction of a triangle-free graph $G$ with $e$ edges and
\begin{equation} bip(G) \leq \frac{e}{2} + \frac{1}{2} e^{5/6}.\end{equation}

In general, the quadrance can be defined in any m-dimensional space $\mathbbm{F}_q^m$ for $m \geq 2$ as follows. 

\begin{Definition}
  The quadrance $Q ( X, Y )$ between the points $X = (x_1, \ldots, x_m)$ and $Y = ( y_1, \ldots ,y_m )$ in $\mathbbm{F}_q^m $ is the number
  
  \[ Q ( X, Y ) := \sum_{i=1}^{m} ( x_i - y_i )^2. \]
\end{Definition}

For any fixed $a \neq 0 \in \mathbbm{F}_q$, the quadrance graph $D_q^{m}(a)$ has the vertex set $\mathbbm{F}_q^m$, and $X, Y \in
\mathbbm{F}_q^m$ are adjacent if and only if $Q ( X, Y ) =a$. Similarly to the above, we have the following theorem.

\begin{Theorem} \label{main theorem 1} Let $q$ be any odd prime power and $a \neq 0 \in \mathbbm{F}_q$. The unit-quadrance graph $D_q^{m}(a)$ is a $d_{q,m}$-regular graph on $n_{q,m}=q^m$ vertices with the following properties. 
\begin{enumerate}
  \item (\cite{before}) Let $\chi$ be the \textit{quadratic character}. Equivalently, $\chi$ is $1$ on squares, $0$ at $0$ and $-
1$ otherwise. Then
  \begin{equation*}
   d_{q,m} =
   \begin{cases}
    q^{m-1} + \chi((-1)^{(m-1)/2})q^{(m-1)/2} & \text{if}\; \; m \;  \text{is odd},\\
    q^{m-1} - \chi((-1)^{n/2})q^{(n-2)/2} & \text{otherwise}.
   \end{cases}
  \end{equation*}
  \item (\cite{before}) Let $\lambda \ne \Delta(D_q^{m})$ be any eigenvalue of the graph $D_q^m$. Then 
         \[ | \lambda | \leq 2 q^{(m-1)/2}. \] 
	\item The toughness of $D_q$ is at least \[ q^{(m-1)/2}/(6+o(1)) = n_{q,m}^{(m-1)/2m}/(6+o(1)).\]
	\item The bisection width of $D_q$ is at least $q^m (d_{q,m}-2q^{(m-1)/2})/(4+o(1))$.
	\item The maximum number of edges in a spanning bipartite subgraph of $D_q$ is at most $bip(D_q) \leq q^m (d_{q,m}+2q^{(m-1)/2}) /4$.
	\item The independence number of $D_q$ is at most $(2+o(1))q^{(m+1)/2} = (2+o(1)) n_{q,m}^{(m+1)/2m}$, and hence its chromatic number is at least $n_q^{(m-1)/2m}/(2+o(1))$ 
\end{enumerate}
\end{Theorem}

The proof of this theorem is omitted since it is the same as the proof of Theorem \ref{main theorem}. Note that $D_q^m(a)$ is triangle-free if and only if $m=2$ and $q$ is a prime of the form $q = 12k \pm 5$. Thus Theorem \ref{main theorem 1} is not more useful than Theorem \ref{main theorem} in attacking triangle-free graphs. 

\section*{\large 5. Properties of finite non-Euclidean graphs}

The following theorem summarizes some of the properties of finite non-Euclidean graphs.
\begin{Theorem} \label{mt1} Let $q$ be any prime of the form $q = 12k + 5$ with $k \geq 1$. The finite non-Euclidean graph $V_q(3,6)$ is a $d_q=(q + 1)$-regular graph on $n_q=q^2-q$ vertices with the following properties. 
\begin{enumerate}
	\item $V_q(3,6)$ is triangle-free and has diameter $3$.
	\item The toughness of $V_q(3,6)$ is at least $q^{1/2}/6 > n_q^{1/4}/6$.
	\item The bisection width of $V_q(3,6)$ is at least $(q^2-q) (q-2q^{1/2})/(4+o(1))$.
	\item The maximum number of edges in a spanning bipartite subgraph of $V_q(3,6)$ is at most $bip(V_q(3,6)) \leq (q^2-q) (q+2q^{1/2}) /4$.
	\item The independence number of $V_q(3,6)$ is at most $(2+o(1))n_q^{3/4}$, and hence its chromatic number is at least $n_q^{1/4}/(2+o(1))$ 
\end{enumerate}
\end{Theorem}

\begin{proof}
Part 1 follows from Theorems \ref{girth1} and \ref{diameter1} combined with the fact that $3$ is not a square in $\mathbbm{F}_q$ if $q = 12k + 5$.
 
Part 2 follows directly from Lemma \ref{eigenvalue 2} and Theorem \ref{tough theorem}. 

Parts 3 and 4 are immediate from Lemma \ref{eigenvalue 2} and Lemma \ref{bip}.

The proof of Part 5 is similar to the proof of Theorem \ref{main theorem}, Part 5.
\end{proof}

Theorem \ref{mt1} shows that the finite non-Euclidean graph $V_q(3,6)$, where $q$ is a prime of the form $q=12k+5$, is an explicit triangle-free graph on $n_q=q^2-q$ vertices whose chromatic number exceeds $(0.5+o(1))n_q^{1/4}$. In addition, the finite non-Euclidean graph $V_q(a)$ is an explicit construction showing that $R(3,k) \geq \Omega(k^{4/3})$.

\section*{\large 6. Further remarks}

It is known that the chromatic number of any graph with maximum degree $d$ in which the
number of edges in the induced subgraph on the set of all neighbors of any vertex does not exceed
$d^2/f$ is at most $O(d / log f)$ (see \cite{alon4}). Let $q$ be a prime of the form $q=12k \pm 5$. Then the induced subgraph on the set of all neighbors of any vertex is empty. So we can set $f = q^2$. 
For any $q$, the induced subgraph on the the set of all neighbors of any vertex of $D_q(a)$ has at most $q+1$ edges. This implies that we can set $f = d-1$. This implies the following upper bound for the chromatic number of quadrance graphs $D_q(a)$

\begin{equation} \chi(D_q(a)) \leq O(q/\log_2q).  \end{equation}

This method also gives us a similar bound for higher dimensional cases.

\begin{equation} \chi(D_q^m(a)) \leq O(q^{m-1}/\log_2q).  \end{equation}

Note that the best known constructive upper bound for the chromatic number of general quadrance graphs is much weaker: 
\begin{equation} \chi(D_q^m) \leq q^{m-1}(1/2 + o(1)). \end{equation}
This bound is obtained by an explicit coloring in \cite{vinh}.

Finally, the bounds in Theorems \ref{main theorem} and \ref{mt1} match the bounds obtained by code graphs in Theorem 3.1 in \cite{alon1}. These latter graphs are Cayley graphs and their construction is based on some of the properties of certain dual BCH error-correcting codes. For a positive integer $k$, let $F_k = GF(2^k)$ denote the finite field with $2^k$ elements. The elements of $F_k$ are represented by binary vectors of length $k$. If $a$ and $b$ are two such vectors, let $(a,b)$ denote their concatenation. Let $G_k$ be the graph whose vertices are all $n = 2^{2k}$ binary vectors of length $2k$, where two vectors $u$ and $v$ are adjacent if and only if there exists a non-zero $z \in F_k$ such that $u+v=(z,z^3)$ mod $2$ where $z^3$ is computed in the field $F_k$. Then $G_k$ is a $d_k=2^k-1$-regular graph on $n_k=2^{2k}$ vertices. Moreover, $G_k$ is triangle-free with independence number at most $2n^{3/4}$. Alon gives the better bound $R(m,3) \geq \Omega(m^{3/2})$ in \cite{alon} by considering a graph whose vertex set is the set of all $n=2^{3k}$ binary vectors of length $3k$ (instead of all binary vectors of length $2k$). Suppose that $k$ is not divisible by $3$. Let $W_0$ be the set of all nonzero elements $\alpha \in F_k$ such that the leftmost bit in the binary representation of $\alpha^7$ is $0$, and let $W_1$ be the set of all nonzero elements $\alpha \in F_k$ for which the leftmost bit of $\alpha^7$ is $1$. Then $|W_0|=2^{k-1}-1$ and $|W_1|=2^{k-1}$. Let $G_n$ be the graph whose vertices are all $n=2^{3k}$ binary vectors of length $3k$, where two vectors $u$ and $v$ are adjacent if and only if there exist $w_0 \in W_0$ and $w_1 \in W_1$ such that $u+v=(w_0,w_0^3,w_0^5)+(w_1,w_1^3,w_1^5)$ where the powers are computed in the field $F_k$ and the addition is addition modulo $2$. Then $G_n$ is a $d_n=2^{k-1}(2^{k-1}-1)$-regular graph on $n=2^{3k}$ vertices. Moreover, $G_n$ is a triangle-free graph with independence number at most $(36+o(1))n^{2/3}$. Note that, going to higher dimensional unit-quadrance graphs does not give us a better bound for the Ramsey number $R(3,k)$ as for any $m \geq 3$ the graph $D_q^m(a)$ is not triangle-free. The problem of finding better bounds for the chromatic number of the graphs $D_q(a)$ and $V_q(\sigma,a)$ touches on an important question in graph theory: what is the greatest possible chromatic number for a triangle-free regular graph of order n? A possible approach is to consider the existence of sum-free varieties in high dimensional vector spaces over finite fields. We see that the varieties of degree two only give us triangle-free graphs on vector spaces of dimension two. We hope to address this problem for varieties of higher degrees in higher dimensions in a subsequent paper.


\nocite{*}
\begin{thebibliography}{1}
\small
\providecommand{\natexlab}[1]{#1}
\providecommand{\url}[1]{\texttt{#1}}
\expandafter\ifx\csname urlstyle\endcsname\relax
  \providecommand{\doi}[1]{doi: #1}\else
  \providecommand{\doi}{doi: \begingroup \urlstyle{rm}\Url}\fi




\bibitem{angel1}
Angel, J., Finite upper half planes over finite fields, \textit{Finite Fields Appl.} \textbf{2} (1996), 62-68.

\bibitem{angel}
Angel, J. and Evans, R., Diameters of finite upper half plane graphs, \textit{Journal of Graph Theory} \textbf{23} (2) (1996), 129-137.

\bibitem{alon}
Alon, N., Explicit Ramsey graphs and orthonormal labellings, \textit{Electronic Journal of Combinatorics} \textbf{1} (1994), R12, 8pp.

\bibitem{alon1}
Alon, N., Tough Ramsey graphs without short cycles, \textit{Journal of Algebraic Combinatorics} \textbf{4} (1995), 189-195.

\bibitem{3}
Alon, N., Krivelevich, M. and Sudakow, B., List coloring of random and pseudo-random graphs, \textit{Combinatorica} \textbf{19} (4) (1999) 453-472.

\bibitem{alon4}
Alon, N., Krivelevich, M. and Sudakov, B., Coloring graphs with sparse neighborhoods, \textit{J. Combinatorial Theory, Ser. B} \textbf{77} (1999), 73-82.

\bibitem{bauer}
Bauer, D., Vandenheuvel, J. and Schmeichel, E., Toughness and triangle-free graphs, \textit{Journal of Combinatorial Theory, Series B}, \textbf{65} (2) (1995) 208-221.

\bibitem{celniker}
Celniker, N., Eigenvalue bounds and girths of graphs of finite, upper half-planes, \textit{Pacific Journal of Mathematics}, \textbf{166} (1), (1994) 1-21.

\bibitem{celniker1}
Celniker, N., Poulos, S., Terras, A., Trimble, C. and Velasquez, E., Is there life on finite upper half planes?, \textit{Contemp. Math.} \textbf{143} (1993), 65-88.

\bibitem{chvatal}
Chv\'atal, V., Tough graphs and hamiltonian circuits, \textit{Discrete Mathematics} \textbf{5} (1973), 215-218.

\bibitem{evans1}
Evans, R., Character sums as orthogonal eigenfunctions of adjacency operators for Cayley graphs, \textit{Contemporary Mathematics}, \textbf{168}, American Mathematical Society, Providence, RI, 1994, 33-50.

\bibitem{evans2}
Evans, R., Spherical functions for finite upper half planes with characteristic 2, \textit{Finite Fields Appl.} \textbf{1} (1995), 376-394.

\bibitem{algebraic graph theory}
Godsil, D. and Royle, G., \textit{Algebraic Graph Theory}, Graduate Texts in Mathematics \textbf{207}, Springer-Verlag New York, 2001.

\bibitem{katz1}
Katz, N. M., Estimates for Soto-Andrade sums, \textit{J. Reine Angew. Math.} \textbf{438} (1993), 143-161.

\bibitem{katz2}
Katz, N. M., A note on exponential sums, \textit{Finite Fields Appl.} \textbf{1} (1995), 395-398.


\bibitem{li}
Li, W.-C.W., \textit{Number Theory with Applications}, World Scientific, River Edge, NJ, 1996.


\bibitem
{before} Medrano, A., Myers, P., Stark, H. M., Terras, A., Finite analogues of Euclidean space,
\textit{Journal of Computational and Applied Mathematics}, \textbf{68} (1996), 221-238.

\bibitem
{moorhouse} Moorhouse, E., \textit{On the chromatic number of planes}, manuscript,\, http://www.uwyo.edu/moorhouse/pub/chromatic.pdf

\bibitem
{old} Problem 6, Problems collected at the AGG2004 Workshop,\, http://www.ipam.ucla.edu/programs/agg2004/

\bibitem{terras1}
Terras, A., Eigenvalue problems related to finite analogues of upper half planes, \textit{Discourses in Math.} \textbf{1} (1991), 237-263.

\bibitem{terras2}
Terras, A., Are finite upper half planes ramanujan?, \textit{DIMACS Series in Discrete Math and Theoritical Computer Science} \textbf{10} (1993), 125-142.

\bibitem
{survey} Terras, A., Survey of spectra of laplacians on finite symmetric spaces, \textit{Experimental Mathematics}, (1996).

\bibitem
{vinh2} Vinh, L. A., Random walk on hypergroup of circles in a finite field, \textit{Proceedings of Australasian Workshop on Combinatorial Algroithms} 2005, University of Ballarat, Australia, 341-351.

\bibitem
{vinh3} Vinh, L. A., Random walk on hypergroup of conics over finite fields, to appear on
\textit{The Global Journal in Pure and Applied Mathematics}.

\bibitem
{vinh} Vinh, L. A., Some coloring problems for unit-quadrance graphs,
\textit{Proceedings of Australian Workshop on Combinatorial Algorithms} 2006, Charles Darwin University, Australia, 361-367.

\bibitem
{vinh1} Vinh, L. A., Quadrance graphs, \textit{Australian Mathematical Society Gazette}, \textbf{33} (5), November 2006, 330-332.

\bibitem{vinh4}
Vinh, L. A., Explicit Ramsey graphs and Erdos distance problem over finite Euclidean and non-Euclidean spaces, \textit{Electronics Journal of Combinatorics}, \textbf{15} (1), 2008, R5, 18 pages.

\bibitem{weil}
Weil, A., On the Riemann hypothesis in function fields, \textit{Proc. Natl. Acad. Sci. USA}, \textbf{27} (1941), 345-347.

\bibitem
{quadrance} Wildberger, N. J.,
\textit{Divine Proportions: Rational Trigonometry to Universal Geometry}, WildEgg, Australia, 2005.
\end{thebibliography}
\end{document}